\documentclass[11pt,twoside,a4paper]{article}
\usepackage{amsmath}
\usepackage{amssymb}
\usepackage[english]{babel}
\renewcommand{\a}{{\bf a}}
\renewcommand{\b}{{\bf b}}
\renewcommand{\c}{{\bf c}}

\newcommand{\A}{{\cal A}}

\newcommand{\D}{{\cal D}}
\renewcommand{\L}{{\cal L}}

\newcommand{\X}{{\bf X}}
\newcommand{\Y}{{\bf Y}}
\newcommand{\Z}{{\bf Z}}
\newcommand{\h}{{\cal H}_{CK}}
\newcommand{\T}{\mathbb{T}_{CK}}
\newcommand{\F}{\mathbb{F}_{CK}}
\renewcommand{\H}{{\cal  H}_{NCK}}
\newcommand{\TT}{\mathbb{T}_{NCK}}
\newcommand{\FF}{\mathbb{F}_{NCK}}
\newcommand{\hfdb}{{\cal H}_{FdB}}
\newcommand{\hfdbD}{{\cal H}_{FdB,D}}

\newcommand{\tun}{\begin{picture}(5,0)(-2,-1)
\put(0,0){\circle*{2}}
\end{picture}}

\newcommand{\tdeux}{\begin{picture}(7,7)(0,-1)
\put(3,0){\circle*{2}}
\put(3,0){\line(0,1){5}}
\put(3,5){\circle*{2}}
\end{picture}}

\newcommand{\ttroisun}{\begin{picture}(15,8)(-5,-1)
\put(3,0){\circle*{2}}
\put(-0.65,0){$\vee$}
\put(6,7){\circle*{2}}
\put(0,7){\circle*{2}}
\end{picture}}
\newcommand{\ttroisdeux}{\begin{picture}(5,12)(-2,-1)
\put(0,0){\circle*{2}}
\put(0,0){\line(0,1){5}}
\put(0,5){\circle*{2}}
\put(0,5){\line(0,1){5}}
\put(0,10){\circle*{2}}
\end{picture}}

\newcommand{\tquatreun}{\begin{picture}(15,12)(-5,-1)
\put(3,0){\circle*{2}}
\put(-0.65,0){$\vee$}
\put(6,7){\circle*{2}}
\put(0,7){\circle*{2}}
\put(3,7){\circle*{2}}
\put(3,0){\line(0,1){7}}
\end{picture}}
\newcommand{\tquatredeux}{\begin{picture}(15,18)(-5,-1)
\put(3,0){\circle*{2}}
\put(-0.65,0){$\vee$}
\put(6,7){\circle*{2}}
\put(0,7){\circle*{2}}
\put(0,14){\circle*{2}}
\put(0,7){\line(0,1){7}}
\end{picture}}
\newcommand{\tquatretrois}{\begin{picture}(15,18)(-5,-1)
\put(3,0){\circle*{2}}
\put(-0.65,0){$\vee$}
\put(6,7){\circle*{2}}
\put(0,7){\circle*{2}}
\put(6,14){\circle*{2}}
\put(6,7){\line(0,1){7}}
\end{picture}}
\newcommand{\tquatrequatre}{\begin{picture}(15,18)(-5,-1)
\put(3,5){\circle*{2}}
\put(-0.65,5){$\vee$}
\put(6,12){\circle*{2}}
\put(0,12){\circle*{2}}
\put(3,0){\circle*{2}}
\put(3,0){\line(0,1){5}}
\end{picture}}
\newcommand{\tquatrecinq}{\begin{picture}(9,19)(-2,-1)
\put(0,0){\circle*{2}}
\put(0,0){\line(0,1){5}}
\put(0,5){\circle*{2}}
\put(0,5){\line(0,1){5}}
\put(0,10){\circle*{2}}
\put(0,10){\line(0,1){5}}
\put(0,15){\circle*{2}}
\end{picture}}

\newcommand{\tcinqun}{\begin{picture}(20,8)(-5,-1)
\put(3,0){\circle*{2}}
\put(-0.5,0){$\vee$}
\put(6,7){\circle*{2}}
\put(0,7){\circle*{2}}
\put(3,0){\line(2,1){10}}
\put(3,0){\line(-2,1){10}}
\put(-7,5){\circle*{2}}
\put(13,5){\circle*{2}}
\end{picture}}
\newcommand{\tcinqdeux}{\begin{picture}(15,14)(-5,-1)
\put(3,0){\circle*{2}}
\put(-0.65,0){$\vee$}
\put(6,7){\circle*{2}}
\put(0,7){\circle*{2}}
\put(3,7){\circle*{2}}
\put(3,0){\line(0,1){7}}
\put(0,7){\line(0,1){7}}
\put(0,14){\circle*{2}}
\end{picture}}
\newcommand{\tcinqtrois}{\begin{picture}(15,15)(-5,-1)
\put(3,0){\circle*{2}}
\put(-0.65,0){$\vee$}
\put(6,7){\circle*{2}}
\put(0,7){\circle*{2}}
\put(3,7){\circle*{2}}
\put(3,0){\line(0,1){7}}
\put(3,7){\line(0,1){7}}
\put(3,14){\circle*{2}}
\end{picture}}
\newcommand{\tcinqquatre}{\begin{picture}(15,14)(-5,-1)
\put(3,0){\circle*{2}}
\put(-0.65,0){$\vee$}
\put(6,7){\circle*{2}}
\put(0,7){\circle*{2}}
\put(3,7){\circle*{2}}
\put(3,0){\line(0,1){7}}
\put(6,7){\line(0,1){7}}
\put(6,14){\circle*{2}}
\end{picture}}
\newcommand{\tcinqcinq}{\begin{picture}(15,19)(-5,-1)
\put(3,0){\circle*{2}}
\put(-0.65,0){$\vee$}
\put(6,7){\circle*{2}}
\put(0,7){\circle*{2}}
\put(6,14){\circle*{2}}
\put(6,7){\line(0,1){7}}
\put(0,14){\circle*{2}}
\put(0,7){\line(0,1){7}}
\end{picture}}
\newcommand{\tcinqsix}{\begin{picture}(15,20)(-7,-1)
\put(3,0){\circle*{2}}
\put(-0.65,0){$\vee$}
\put(6,7){\circle*{2}}
\put(0,7){\circle*{2}}
\put(-3.65,7){$\vee$}
\put(3,14){\circle*{2}}
\put(-3,14){\circle*{2}}
\end{picture}}
\newcommand{\tcinqsept}{\begin{picture}(15,8)(-5,-1)
\put(3,0){\circle*{2}}
\put(-0.65,0){$\vee$}
\put(6,7){\circle*{2}}
\put(0,7){\circle*{2}}
\put(2.35,7){$\vee$}
\put(3,14){\circle*{2}}
\put(9,14){\circle*{2}}
\end{picture}}
\newcommand{\tcinqhuit}{\begin{picture}(15,26)(-5,-1)
\put(3,0){\circle*{2}}
\put(-0.65,0){$\vee$}
\put(6,7){\circle*{2}}
\put(0,7){\circle*{2}}
\put(0,14){\circle*{2}}
\put(0,7){\line(0,1){7}}
\put(0,21){\circle*{2}}
\put(0,14){\line(0,1){7}}
\end{picture}}
\newcommand{\tcinqneuf}{\begin{picture}(15,26)(-5,-1)
\put(3,0){\circle*{2}}
\put(-0.65,0){$\vee$}
\put(6,7){\circle*{2}}
\put(0,7){\circle*{2}}
\put(6,14){\circle*{2}}
\put(6,7){\line(0,1){7}}
\put(6,21){\circle*{2}}
\put(6,14){\line(0,1){7}}
\end{picture}}
\newcommand{\tcinqdix}{\begin{picture}(15,19)(-5,-1)
\put(3,5){\circle*{2}}
\put(-0.5,5){$\vee$}
\put(6,12){\circle*{2}}
\put(0,12){\circle*{2}}
\put(3,0){\circle*{2}}
\put(3,0){\line(0,1){12}}
\put(3,12){\circle*{2}}
\end{picture}}
\newcommand{\tcinqonze}{\begin{picture}(15,26)(-5,-1)
\put(3,5){\circle*{2}}
\put(-0.65,5){$\vee$}
\put(6,12){\circle*{2}}
\put(0,12){\circle*{2}}
\put(3,0){\circle*{2}}
\put(3,0){\line(0,1){5}}
\put(0,12){\line(0,1){7}}
\put(0,19){\circle*{2}}
\end{picture}}
\newcommand{\tcinqdouze}{\begin{picture}(15,26)(-5,-1)
\put(3,5){\circle*{2}}
\put(-0.65,5){$\vee$}
\put(6,12){\circle*{2}}
\put(0,12){\circle*{2}}
\put(3,0){\circle*{2}}
\put(3,0){\line(0,1){5}}
\put(6,12){\line(0,1){7}}
\put(6,19){\circle*{2}}
\end{picture}}
\newcommand{\tcinqtreize}{\begin{picture}(5,26)(-2,-1)
\put(0,0){\circle*{2}}
\put(0,0){\line(0,1){7}}
\put(0,7){\circle*{2}}
\put(0,7){\line(0,1){7}}
\put(0,14){\circle*{2}}
\put(-3.65,14){$\vee$}
\put(-3,21){\circle*{2}}
\put(3,21){\circle*{2}}
\end{picture}}
\newcommand{\tcinqquatorze}{\begin{picture}(9,26)(-5,-1)
\put(0,0){\circle*{2}}
\put(0,0){\line(0,1){5}}
\put(0,5){\circle*{2}}
\put(0,5){\line(0,1){5}}
\put(0,10){\circle*{2}}
\put(0,10){\line(0,1){5}}
\put(0,15){\circle*{2}}
\put(0,15){\line(0,1){5}}
\put(0,20){\circle*{2}}
\end{picture}}

\newcommand{\tdun}[1]{\begin{picture}(10,5)(-2,-1)
\put(0,0){\circle*{2}}
\put(3,-2){\tiny #1}
\end{picture}}

\newcommand{\tddeux}[2]{\begin{picture}(12,5)(0,-1)
\put(3,0){\circle*{2}}
\put(3,0){\line(0,1){5}}
\put(3,5){\circle*{2}}
\put(6,-2){\tiny #1}
\put(6,3){\tiny #2}
\end{picture}}

\newcommand{\tdtroisun}[3]{\begin{picture}(20,12)(-5,-1)
\put(3,0){\circle*{2}}
\put(-0.65,0){$\vee$}
\put(6,7){\circle*{2}}
\put(0,7){\circle*{2}}
\put(5,-2){\tiny #1}
\put(9,5){\tiny #2}
\put(-5,5){\tiny #3}
\end{picture}}
\newcommand{\tdtroisdeux}[3]{\begin{picture}(12,12)(-2,-1)
\put(0,0){\circle*{2}}
\put(0,0){\line(0,1){5}}
\put(0,5){\circle*{2}}
\put(0,5){\line(0,1){5}}
\put(0,10){\circle*{2}}
\put(3,-2){\tiny #1}
\put(3,3){\tiny #2}
\put(3,9){\tiny #3}
\end{picture}}

\title{Fa\`a di Bruno subalgebras of the Hopf algebra of planar trees
from combinatorial Dyson-Schwinger equations}
\date{}
\author{Lo{\"\i}c Foissy \\
\\
{\small{\it Laboratoire de Math\'ematiques - UMR6056, Universit\'e de Reims}}\\
\small{{\it Moulin de la Housse - BP 1039 - 51687 REIMS Cedex 2, France}}\\
\small{e-mail : loic.foissy@univ-reims.fr}}

\newtheorem{theo}{\indent Theorem}
\newtheorem{lemme}[theo]{\indent Lemma}
\newtheorem{prop}[theo]{\indent Proposition}
\newtheorem{cor}[theo]{\indent Corollary}
\newtheorem{defi}[theo]{\indent Definition}

\begin{document}

\maketitle

ABSTRACT. We consider the combinatorial Dyson-Schwinger
equation $X=B^+(P(X))$ in the non-commutative Connes-Kreimer
Hopf algebra of planar rooted trees $\H$, where $B^+$ is the operator 
of grafting on a root, and $P$ a formal series.
The unique solution $X$ of this equation generates a graded subalgebra $\A_{N,P}$ of $\H$.
We describe all the formal series $P$ such that $\A_{N,P}$ is a Hopf subalgebra.
We obtain in this way a $2$-parameters family of Hopf subalgebras of $\H$,
organized into three isomorphism classes: 
a first one, restricted to a polynomial ring in one variable; a second one,
restricted to the Hopf subalgebra of ladders, isomorphic to the Hopf algebra
of quasi-symmetric functions; a last (infinite) one, which gives a non-commutative
version of the Fa\`a di Bruno Hopf algebra. By taking the quotient, 
the last classe gives an infinite set of embeddings of the Fa\`a di Bruno
algebra into the Connes-Kreimer Hopf algebra of rooted trees.
Moreover, we give an embedding of the free Fa\`a di Bruno Hopf algebra
on $D$ variables into a Hopf algebra of decorated rooted trees, together
with a non commutative version of this embedding.

\tableofcontents

\section*{Introduction}

The Connes-Kreimer Hopf algebra $\h$ of rooted trees is introduced in \cite{Kreimer1}.
It is commutative and not cocommutative.
A particular Hopf subalgebra of $\h$, namely the Connes-Moscovici subalgebra, is 
introduced in \cite{Connes}. It is the subalgebra generated by the following elements:
$$  \left\{ \begin{array}{rcl}
\delta_1&=&\tun,\\
\delta_2&=&\tdeux,\\
\delta_3&=&\ttroisun+\ttroisdeux,\\
\delta_4&=&\tquatreun+3\tquatredeux+\tquatrequatre+\tquatrecinq,\\
\delta_5&=&\tcinqun+6\tcinqdeux+3\tcinqcinq+4\tcinqsix+4\tcinqhuit+\tcinqdix
+3\tcinqonze+\tcinqtreize+\tcinqquatorze,\\
\vdots&&
\end{array}
\right.$$
The appearing coefficients, called Connes-Moscovici coefficients, are studied in \cite{Brouder,Hoffman}. 
It is shown in \cite{Moscovici} that the character group of this subalgebra is isomorphic to 
the group of formal diffeomorphisms, that is to say the group of formal series of the form
$h+a_1h^2+\ldots$, with composition. In other terms, the Connes-Moscovici subalgebra is isomorphic
to the Hopf algebra of functions on the group of formal diffeomorphisms, also called
the Fa\`a di Bruno Hopf algebra.

A non commutative version $\H$ of the Connes-Kreimer Hopf algebra of trees is introduced
in \cite{Foissy2,Holtkamp}. It contains a non commutative version of the Connes-Moscovici subalgebra,
described in \cite{Foissy3}. Its abelianization can be identified with the subalgebra of $\h$,
here denoted by $\A_{1,1}$, generated by the following elements of $\h$:
$$ \left\{ \begin{array}{rcl}
a_1&=&\tun,\\
a_2&=&\tdeux,\\
a_3&=&\ttroisun+\ttroisdeux,\\
a_4&=&\tquatreun+2\tquatredeux+\tquatrequatre+\tquatrecinq,\\
a_5&=&\tcinqun+3\tcinqdeux+\tcinqcinq+2\tcinqsix+2\tcinqhuit+\tcinqdix
+2\tcinqonze+\tcinqtreize+\tcinqquatorze,\\
\vdots&&
\end{array}
\right.$$
This subalgebra is different from the Connes-Moscovici subalgebra, but is also
isomorphic to the Fa\`a di Bruno Hopf algebra.\\

In this paper, we consider a family of subalgebras of $\H$, which give a non commutative version of the Fa\`a di Bruno algebra. 
They are generated by a combinatorial Dyson-Schwinger equation \cite{Bergbauer,Kreimer5,Yeats}:
$$\X_P=B^+(P(\X_P)),$$
where $B^+$ is the operator of grafting on a common root, and $P=\sum p_k h^k$ is a formal series
such that $p_0=1$.  All this makes sense in a completion of $\H$, where this equation admits a unique solution  $\X_P=\sum \a_k$,
whose coefficients are inductively defined by:
$$  \left\{ \begin{array}{rcl}
\a_1&=&\tun,\\
\a_{n+1}&=&\displaystyle \sum_{k=1}^{n}\: \sum_{\alpha_1+\ldots+\alpha_k=n} 
p_k B^+(\a_{\alpha_1}\ldots \a_{\alpha_k}),
\end{array}
\right.$$
For the usual Dyson-Schwinger equation, $P=\alpha (1-h)^{-1}$.
We characterise the formal series $P$ such that the associated subalgebra is Hopf:
we obtain a two-parameters family $\A_{N,\alpha,\beta}$ of Hopf subalgebras of $\H$
and we explicitely describe the system of generator of these algebras.

We then characterise the equalities between the $\A_{N,\alpha,\beta}$'s and then their isomorphism classes.
We obtain three classes:
\begin{enumerate}
\item $\A_{N,0,1}$, equal to $K[\tun]$.
\item $\A_{N,1,-1}$, the subalgebra of ladders, isomorphic to the Hopf algebra of quasi-symmetric functions.
\item The $\A_{N,1,\beta}$'s,  with $\beta\neq -1$, a non commutative version of the Fa\`a di Bruno Hopf algebra.
\end{enumerate}
By taking the quotient, we obtain three classes of Hopf subalgebras of $\h$:
\begin{enumerate}
\item $\A_{0,1}$, equal to $K[\tun]$.
\item $\A_{1,-1}$, the subalgebra of ladders, isomorphic to the Hopf algebra of symmetric functions.
\item The $\A_{1,\beta}$'s,  with $\beta\neq -1$, isomorphic to the Fa\`a di Bruno Hopf algebra.
\end{enumerate}

We finally give an embedding of a non commutative version of the free Fa\`a di Bruno on $D$ variables (see \cite{Effros})
in a Hopf algebra of planar rooted trees decorated by the set $\{1,\ldots, D\}^3$. By taking the quotient,
the free Fa\`a di Bruno algebra appears as a subalgebra of a Hopf algebra of decorated rooted trees.\\

This text is organized as follows. The first section gives some recalls about the Hopf algebras of trees
and the Fa\`a di Bruno algebra. We define the subalgebras of $\h$ and $\H$ associated to a formal series
$P$ in section 2 and also give here the main theorem (theorem \ref{theoprincipal}),
which characterizes the $P$'s such that the associated subalgebras are Hopf.
In section 3, we prove $2\Longrightarrow 3$ of theorem \ref{theoprincipal}.
In section 4, we prove $4\Longrightarrow 1$ of theorem \ref{theoprincipal}.
We also describe there the system of generators, and the case of equalities of the subalgebras.
We describe the isomorphism classes of these subalgebras in the following section.
In the last one, we consider the multivariable case.\\

{\bf Notations.}
\begin{enumerate}
\item $K$ is any field of characteristic zero.
\item Let $\lambda \in K$. We put:
$$g_\lambda(h)=(1-h)^{-\lambda}
=\sum_{k=0}^\infty \frac{\lambda(\lambda+1)\ldots (\lambda+k-1)}{k!} h^k 
=\sum_{k=0}^\infty Q_k(\lambda) h^k \in K[[h]].$$
\end{enumerate}

\section{Preliminaries}

\subsection{Valuation and $n$-adic topology}

In this paragraph, let us consider a graded Hopf algebra $\A$.
Let $\A_n$ be the homogeneous component of degree $n$ of $\A$.
For all $a \in \A$, we put:
$$val(a)=\max\left\{n \in \mathbb{N}\:/\: a \in \bigoplus_{k\geq n} \A_k\right\}
\in \mathbb{N}\cup\{+\infty\}.$$
For all $a,b \in \A$, we also put $d(a,b)=2^{-val(a-b)}$,
with the convention $2^{-\infty}=0$. 
Then $d$ is a distance on $\A$.
The induced topology over $\A$ will be called the {\it $n$-adic topology}.\\

Let $\overline{\A}$ be the completion of $\A$ for this distance.
In other terms:
$$\overline{\A}=\prod_{n=0}^{+\infty} \A_n.$$
The elements of  $\overline{\A}$ are written in the form 
$\displaystyle \sum_{n=0}^{+\infty} a_n$, with $a_n \in \A_n$ for all $n$.
Moreover, $\overline{\A}$ is naturally given a structure of associative algebra,
by continuously extending the product of $\A$. 
The coproduct of $\A$ can also be extended in the following way:
$$\Delta : \overline{\A}  \longrightarrow  \overline{\A} \hat{\otimes}  \overline{\A}
=\prod_{i,j \in \mathbb{N}} \A_i \otimes \A_j. $$

For all $\displaystyle p=\sum_{k=0}^{+\infty} a_n h_n \in K[[h]]$, for all
$a \in \overline{\A}$ such that $val(a)\geq 1$ , we put:
$$p(a)=\sum_{k=0}^{+\infty} p_n a^n \in \overline{\A}.$$
Indeed, for all $n,m\in \mathbb{N}$, $\displaystyle val\left(\sum_{k=n}^{n+m}p_n a^n \right)\geq n$,
so this series is Cauchy, and converges.
It is an easy exercise to prove that for all $p,q \in K[[h]]$, such that $q$ has no constant term,
for all $a\in \overline{\A}$, with $val(a)\geq 1$,
$(p\circ q)(a)=p(q(a))$.

\subsection{The commutative Connes-Kreimer Hopf algebra of trees}

This Hopf algebra is introduced by Kreimer in 
\cite{Kreimer1} and studied for example in \cite{Kreimer,Connes,Hoffman,Foissy,Kreimer2,Kreimer3}.

\begin{defi}
\textnormal{\begin{enumerate}
\item A {\it rooted tree} $t$ is a finite, connected graph without loops,
with a special vertex called {\it root}.
The set of rooted trees will be denoted by $\T$.
\item The {\it weight} of a rooted tree is the number of its vertices.
\item A {\it planar rooted tree} is a tree which is given an imbedding in the plane.
The set of planar rooted trees will be denoted by $\TT$.
\end{enumerate}}
\end{defi}

{\bf Examples.} \begin{enumerate}
\item Rooted trees of weight $\leq 5$:
$$\tun,\tdeux,\ttroisun,\ttroisdeux,\tquatreun, \tquatredeux,
\tquatrequatre,\tquatrecinq,
\tcinqun,\tcinqdeux,\tcinqcinq,\tcinqsix,\tcinqhuit,
\tcinqdix,\tcinqonze,\tcinqtreize,\tcinqquatorze.$$
\item Planar rooted trees of weight $\leq 5$:
$$\tun,\tdeux,\ttroisun,\ttroisdeux,\tquatreun, \tquatredeux,
\tquatretrois,\tquatrequatre,\tquatrecinq,
\tcinqun,\tcinqdeux,\tcinqtrois,\tcinqquatre,\tcinqcinq,\tcinqsix,\tcinqsept,\tcinqhuit,
\tcinqneuf,\tcinqdix,\tcinqonze,\tcinqdouze,\tcinqtreize,\tcinqquatorze.$$
\end{enumerate}

The Connes-Kreimer Hopf algebra of rooted trees $\h$ is the free associative commutative
algebra freely generated over $K$ by the elements of $\T$.
A linear basis of $\h$ is given by {\it rooted forests}, that is to say monomials in
rooted trees. 
The set of rooted forests will be denoted by $\F$. The weight of a rooted forest
$F=t_1\ldots t_n$ is the sum of the weights of the $t_i$'s.\\

{\bf Examples.} Rooted forests of weight $\leq 4$:
$$1,\tun,\tun\tun,\tdeux,\tun\tun\tun,\tdeux\tun,\ttroisun,\ttroisdeux,
\tun\tun\tun\tun,\tdeux\tun\tun,\ttroisun\tun,\ttroisdeux\tun,\tdeux\tdeux,
\tquatreun,\tquatredeux,\tquatrequatre,\tquatrecinq.$$

We now recall the Hopf algebra structure of $\h$.
An {\it admissible cut} of $t$ is a non empty cut
such that every path in the tree meets at most one cut edge.
The set of admissible cuts of $t$ is denoted by ${\cal A}dm(t)$.
If $c$ is an admissible cut of $t$, 
one of the trees obtained after the application of $c$ contains the root of $t$:
we shall denote it by $R^c(t)$. The product of the other trees will be denoted by $P^c(t)$.
The coproduct of $t$ is then given by :
$$\Delta(t)=t\otimes 1+1\otimes t+\sum_{c \in {\cal A}dm(t)} P^c(t)\otimes R^c(t).$$
This coproduct is extended as an algebra morphism. Then $\h$ becomes a Hopf algebra.
Note that $\h$ is given a gradation of Hopf algebra by the weight.\\

{\bf Examples.} 
\begin{eqnarray*}
\Delta(\tquatreun)&=&\tquatreun \otimes 1+1 \otimes \tquatreun
+ 3\tun\otimes \ttroisun+3 \tun\tun \otimes \tdeux+\tun\tun\tun \otimes \tun,\\
\Delta(\tquatredeux)&=&\tquatredeux \otimes 1+1 \otimes \tquatredeux
+\tdeux\tun\otimes \tun+\tdeux\otimes \tdeux
+ \tun \otimes \ttroisdeux+\tun\tun\otimes \tdeux +\tun \otimes \ttroisun,\\
\Delta(\tquatrequatre)&=&\tquatrequatre\otimes1+1 \otimes \tquatrequatre
+\ttroisun \otimes \tun+ \tun\tun \otimes \tdeux+2 \tun \otimes \ttroisdeux,\\
\Delta(\tquatrecinq)&=&\tquatrecinq \otimes1+1 \otimes \tquatrecinq
+\ttroisdeux\otimes \tun+\tdeux\otimes \tdeux+\tun \otimes \ttroisdeux.
\end{eqnarray*}

We define the operator $B^+:\h\longrightarrow \h$,
that associates to a forest $F\in \F$ the tree obtained by grafting
the roots of the trees of $F$ on a common root.
For example, $B^+(\tdeux\tun)=\tquatredeux$.
Then, for all $x\in \h$:
\begin{equation}
\label{E1}\Delta(B^+(x))=B^+(x)\otimes 1+(Id\otimes B^+)\circ \Delta(x).
\end{equation}
This means that $B^+$ is a $1$-cocycle for a certain cohomology of coalgebra, 
see \cite{Connes} for more details.
Moreover, this operator $B^+$ is homogeneous of degree $1$,
so is continuous. So it can be extended in an operator
$B^+: \overline{\h} \longrightarrow \overline{\h}$.

\subsection{The non commutative Connes-Kreimer Hopf algebra of planar trees}

This algebra is introduced simultaneously in \cite{Foissy2,Holtkamp}.
As an algebra, $\H$ is the free associative algebra generated by the elements of $\TT$.
A basis of $\H$ is given by planar rooted forests, that is to say words in elements of $\TT$.
The set of planar rooted forests will be denoted by $\FF$.\\

{\bf Examples.} Planar rooted forests of weight $\leq 4$:
$$1,\tun,\tun\tun,\tdeux,\tun\tun\tun,\tdeux\tun,\tun \tdeux,\ttroisun,\ttroisdeux,
\tun\tun\tun\tun,\tdeux\tun\tun,\tun \tdeux \tun, \tun \tun \tdeux,
\ttroisun\tun,\tun \ttroisun,\ttroisdeux\tun,\tun \ttroisdeux,\tdeux\tdeux,
\tquatreun,\tquatredeux,\tquatretrois,\tquatrequatre,\tquatrecinq.$$

The coproduct of $\H$ is defined, as for $\h$, with the help
of admissible cuts. For example :

\begin{eqnarray*}
\Delta(\tquatredeux)&=&\tquatredeux \otimes 1+1 \otimes \tquatredeux
+\tdeux\tun\otimes \tun+\tdeux\otimes \tdeux
+ \tun \otimes \ttroisdeux+\tun\tun\otimes \tdeux +\tun \otimes \ttroisun,\\
\Delta(\tquatretrois)&=&\tquatretrois \otimes 1+1 \otimes \tquatretrois
+\tun\tdeux\otimes \tun+\tdeux\otimes \tdeux
+ \tun \otimes \ttroisdeux+\tun\tun\otimes \tdeux +\tun \otimes \ttroisun.
\end{eqnarray*}
Note that $\H$ is a graded Hopf algebra, with a gradation given by the weight.\\

We define an operator, also denoted by $B^+:\H\longrightarrow \H$, as for $\h$. 
For example, $B^+(\tun \tdeux)=\tquatretrois$ and $B^+(\tdeux \tun)=\tquatredeux$.
Then (\ref{E1}) is also satisfied on $\H$.
Moreover, this operator $B^+$ is homogeneous of degree $1$,
so is continuous. In consequence, it can be extended in an operator
$B^+: \overline{\H} \longrightarrow \overline{\H}$.\\

We proved in \cite{Foissy2,Foissy3} that $\H$ is a self-dual Hopf algebra:
it has a non degenerate pairing denoted by $<,>$, and a dual basis $(e_F)_{F \in \FF}$
of the basis of planar decorated forests. The product in the dual basis is given by graftings
(in an extended sense).
For example: 
$$e\tdeux . e \tun\tdeux =e\tdeux\tun \tdeux+e \ttroisdeux \tdeux+
e\tun \tdeux \tdeux+e \tun \tquatredeux+e \tun \tquatrecinq+e\tun \tquatretrois+e\tun \tdeux \tdeux.$$

\subsection{The Fa\`a di Bruno Hopf algebra}

Let $K[[h]]$ be the ring of formal series in one variable over $K$.
We consider:
$$G=\left\{h+\sum_{n\geq 1} a_n h^{n+1} \in K[[h]]\right\}.$$
This is a group for the composition of formal series.
The Fa\`a di Bruno Hopf algebra $\hfdb$ is the Hopf algebra of functions on the opposite of the group $G$.
More precisely, $\hfdb$ is the polynomial ring in variables $Y_i$, with $i\in \mathbb{N}^*$, 
where $Y_i$ is the function on $G$ defined by:
$$Y_i: \left\{ \begin{array}{rcl}
G&\longrightarrow & K\\
\displaystyle h+\sum_{n\geq 1} a_n h^{n+1}&\longrightarrow & a_i.
\end{array}
\right.$$
The coproduct is defined in the following way: for all $f \in \hfdb$, for all $P,Q \in G$,
$$\Delta(f)(P\otimes Q)=f(Q\circ P).$$

This Hopf algebra is commutative and not cocommutative. It is also a graded, connected Hopf algebra,
with $Y_i$ homogeneous of degree $i$ for all $i$. 
We put:
$$\Y=1+\sum_{n=1}^{\infty} Y_n \in \overline{\hfdb}.$$
Then:
$$\Delta(\Y)=\sum_{n=1}^{\infty} \Y^{n+1}\otimes Y_n.$$
Indeed, with the convention $a_0=b_0=1$:
\begin{eqnarray*}
\Delta(\Y)\left(\left(\sum_{i\geq 0} a_i h^{i+1}\right)
\otimes \left(\sum_{i\geq 0} b_i h^{i+1}\right) \right)
&=&\Y \left(\sum_{i\geq 0} b_i \left(\sum_{j\geq 0} a_j h^{j+1}\right)^{i+1}\right)\\
&=&\sum_{i\geq 0} b_i \left(\sum_{j\geq 0} a_j \right)^{i+1},\\
\\
\left(\sum_{n=1}^{\infty} \Y^{n+1}\otimes Y_n\right)\left(\left(\sum_{i\geq 0} a_i h^{i+1}\right)
\otimes \left(\sum_{i\geq 0} b_i h^{i+1}\right) \right)
&=&\sum_{n=1}^{\infty}\left(\sum_{i\geq 0} a_i \right)^{n+1}b_n.\\
\end{eqnarray*}

The graded dual $\hfdb^*$ is an enveloping algebra,
by the Cartier-Quillen-Milnor-Moore theorem. A basis of $Prim(\hfdb^*)$ is given by $(Z_i)_{i\in \mathbb{N}^*}$,
where:
$$ Z_i: \left\{ \begin{array}{rcl}
\hfdb&\longrightarrow & K\\
Y_1^{\alpha_1}\ldots Y_k^{\alpha_k}&\longrightarrow &0\mbox{ if } \alpha_1+\ldots+\alpha_k\neq 1,\\
Y_j&\longrightarrow& \delta_{i,j}.
\end{array}
\right.$$
By homogeneity, for all $i,j \in \mathbb{N}^*$,
there exists a coefficient $\lambda_{i,j} \in K$ such that $[Z_i,Z_j]=\lambda_{i,j} Z_{i+j}$.
Moreover:
\begin{eqnarray*}
[Z_i,Z_j](\Y)&=&\lambda_{i,j}\\
&=&(Z_i\otimes Z_j-Z_j \otimes Z_i)\circ\Delta(\Y)\\
&=&(Z_i\otimes Z_j-Z_j \otimes Z_i)\left(\sum_{n=1}^{\infty} \Y^{n+1}\otimes Y_n\right)\\
&=&Z_i(\Y^{j+1})-Z_j(\Y^{i+1})\\
&=&(j+1)-(i+1)\\
&=&j-i.
\end{eqnarray*}
So the bracket of $Prim(\hfdb^*)$ is given by:
$$[Z_i,Z_j]=(j-i)Z_{i+j}.$$

\section{Subalgebras associated to a formal series}

\subsection{Construction}

We denote by $K[[h]]_1$ the set of formal series of $K[[h]]$ with constant term equal to $1$.

\begin{prop}
Let $P \in K[[h]]_1$.
\begin{enumerate}
\item  There exists a unique element $X_P=\displaystyle \sum_{k\geq 1}a_n \in \overline{\h}$, 
such that $X_P=B^+\left(P(X_P)\right)$.
\item There exists a unique element $\X_P=\displaystyle \sum_{k\geq 1}a_n \in \overline{\H}$,
such that $\X_P=B^+\left(P(\X_P)\right)$.
\end{enumerate}
\end{prop}

{\bf Proof.} 
\begin{enumerate}
\item {\it Unicity.} We put $\displaystyle X_P=\sum_{n \geq 1}a_n$, with $a_n$
homogeneous of degree $n$ for all $n$. Then the $a_n$'s satisfy the following equations:
\begin{eqnarray}
\label{eq1}
&&\left\{ \begin{array}{rcl}
a_1&=&\tun,\\
a_{n+1}&=&\displaystyle \sum_{k=1}^{n}\: \sum_{\alpha_1+\ldots+\alpha_k=n} 
p_k B^+(a_{\alpha_1}\ldots a_{\alpha_k}).
\end{array}\right.
\end{eqnarray}
Hence, the $a_n$'s are uniquely defined. 

{\it Existence.} The $a_n$'s defined inductively by (\ref{eq1}) satisfy the required condition.

\item We put $\displaystyle \X_P=\sum_{n \geq 1}\a_n$, $\a_n$
homogeneous of degree $n$ for all $n$. Then the $\a_n$'s satisfy the following equations:
\begin{eqnarray}
\label{eq1bis}
\left\{ \begin{array}{rcl}
\a_1&=&\tun,\\
\a_{n+1}&=&\displaystyle \sum_{k=1}^{n}\: \sum_{\alpha_1+\ldots+\alpha_k=n} 
p_k B^+(\a_{\alpha_1}\ldots \a_{\alpha_k}).
\end{array}
\right.
\end{eqnarray}
The end of the proof is similar.  $\Box$
\end{enumerate}

\begin{defi}
\textnormal{Let $P \in K[[h]]_1$. 
\begin{enumerate}
\item The subalgebra $\A_P$ of $\h$ is the subalgebra generated by the $a_n$'s.
\item The subalgebra $\A_{N,P}$ of $\H$ is the subalgebra generated by the $\a_n$'s.
\end{enumerate}}
\end{defi}

{\bf Remarks.}
\begin{enumerate}
\item  $\A_P$ is a graded subalgebra of $\h$,
and $\A_{N,P}$ is a graded subalgebra of $\H$.
\item For all $n \in \mathbb{N}^*$, $a_n$ is an element of $Vect(\T)$.
Hence:
\begin{equation}
\label{eq2}
Vect(\T) \cap \A_P=Vect(a_n,\:n\in \mathbb{N}^*).
\end{equation}
The same holds in the non commutative case.
\end{enumerate}

\subsection{Main theorem}

One of the aim of this paper is to prove the following theorem:

\begin{theo}
\label{theoprincipal}
Let $P \in K[[h]]_1$. The following assertions are equivalent:
\begin{enumerate}
\item $\A_{N,P}$ is a Hopf subalgebra of $\H$.
\item $\A_P$ is a Hopf subalgebra of $\h$.
\item There exists $(\alpha,\beta) \in K^2$, 
such that $P$ satisfies the following differential system:
$$ S_{\alpha,\beta}: \left\{\begin{array}{rcl}
(1-\alpha\beta h)P'(h)&=&\alpha P(h)\\
P(0)&=&1.
\end{array}\right. $$
\item There exists $(\alpha,\beta) \in K^2$, such that:
\begin{enumerate}
\item $P(h)=1$ if $\alpha=0$.
\item $P(h)=e^{\alpha h}$ if $\beta=0$.
\item $\displaystyle P(h)=(1-\alpha\beta h)^{-\frac{1}{\beta}}$ if $\alpha \beta \neq 0$.
\end{enumerate}
\end{enumerate}
\end{theo}

An easy computation proves the equivalence between assertions $3$ and $4$.
Moreover, using the Hopf algebra morphism $\varpi:\H\longrightarrow \h$, 
defined by forgetting the planar data, it is clear that $\A_P=\varpi(\A_{N,P})$.
So, assertion $1$ implies assertion $2$.

\section{When is $\A_P$ a Hopf subalgebra? }

\subsection{Preliminary results}
The aim of this section is to show $2\Longrightarrow 3$ in 
theorem \ref{theoprincipal}.

\begin{lemme}
\label{lemme29}
Suppose that $\A_P$ is a Hopf subalgebra. Then two cases are possible:
\begin{enumerate}
\item $P=1$. In this case,  $X_P=\tun$ and $\A_P=K[\tun].$
\item $p_1 \neq 0$. In this case, $a_n\neq 0$ for all $n\geq 1$.
\end{enumerate}
\end{lemme}

{\bf Proof.} Suppose first that $p_1=0$,
and suppose that there exists $n\geq 2$ such that $p_n \neq 0$.
Let us choose $n$ minimal. Then, by (\ref{eq1}), 
$a_2=\ldots=a_n=0$ and $a_{n+1}=p_n B^+(\tun^n)$.
Then:
$$\Delta(B^+(\tun^n))=B^+(\tun^n)\otimes 1+\sum_{k=0}^n \binom{n}{k}
\tun^k\otimes B^+(\tun^{n-k})
\in  \A_P\otimes \A_P,$$
which implies that $B^+(\tun) =\tdeux \in \A_P\cap Vect(\T)=Vect(a_n)$, 
so $a_2 \neq 0$: contradiction. So $P=1$ and $X_P=\tun$.

Suppose $p_1\neq 0$. By (\ref{eq1}), the canonical projection of $a_{n+1}$ on 
$Im((B^+)^2)$ (vector space of trees such that the root has only one child)
is $p_1 B^+(a_n)$ for all $n \geq 1$. Hence, for all $n \geq 1$, $a_{n+1}=0 \Rightarrow a_n=0$.
So for all $n \geq 1$, $a_n \neq 0$. $\Box$\\

We put:
$$ Z: \left\{\begin{array}{rcl}
\h & \longrightarrow &K\\
F\in \F & \longrightarrow & \delta_{\tun,F}.
\end{array}\right. $$

Note that $Z$ is an element of the graded dual $\h^*$.
Moreover, $Z$ can be extended to $\overline{\h}$, and satisfies, 
for all $a,b \in \overline{\h}$:
$$Z(ab)=Z(a) \varepsilon(b) +\varepsilon(a)Z(b).$$

\begin{lemme}
\label{lemme27}
Let $P \in K[[h]]_1$. If $\A_P$ is a Hopf subalgebra of $\h$, then:
$$(Z \otimes Id) \circ \Delta(X_P) \in \overline{\A_P}.$$
\end{lemme}

{\bf Proof.} As $\A_P$ is a Hopf subalgebra, 
for all $n \in \mathbb{N}^*$, $\Delta(a_n) \in \A_P\otimes \A_P$.
Hence, $(Z \otimes Id)\circ\Delta(a_n) \in \A_P$. $\Box$

\begin{lemme}
We consider the following continuous applications:
$$\left\{ \begin{array}{rcccccl}
\tilde{Z}&=&Z \hat{\otimes} Id&:& \overline{\h} 
\hat{\otimes} \overline{\h} &\longrightarrow &\overline{\h},\\
\tilde{\varepsilon}&=&\varepsilon \hat{\otimes} Id&:& \overline{\h} 
\hat{\otimes} \overline{\h} &\longrightarrow& \overline{\h}.
\end{array}
\right.$$

Then $\tilde{\varepsilon}$ is an algebra morphism 
and $\tilde{Z}$ is a $\tilde{\varepsilon}$-derivation, i.e. satisfies:
$$\tilde{Z}(ab)= \tilde{Z}(a) \tilde{\varepsilon}(b) +\tilde{\varepsilon}(a)\tilde{Z}(b) .$$
\end{lemme}

{\bf Proof.} Immediate. $\Box$\\

Let us fix $t \in \T$. We put $\displaystyle P(h)=\sum_{n=0}^{+\infty}p_n h^n.$
As $X_P=B^+(P(X_P))$, we have:
\begin{eqnarray*}
\tilde{Z} \circ \Delta(X_P) &=&\sum_{n=0}^{+\infty} p_n 
(Z \otimes Id) \circ \Delta \circ B^+(X_P^n)\\
 &=&\sum_{n=0}^{+\infty} p_n  Z(B^+(X_P^n)) 1 
+\sum_{n=0}^{+\infty} p_n 
(Z \otimes B^+) (\Delta (X_P))^n\\
 &=&Z(X_P)1 
+B^+\left(\sum_{n=0}^{+\infty} p_n 
\tilde{Z}(\Delta(X_P)^n)\right)\\
 &=&Z(X_P) 1
+B^+\left(\sum_{n=0}^{+\infty} np_n \tilde{\varepsilon}(\Delta(X_P))^{n-1}
\tilde{Z}(\Delta(X_P))
\right)\\
 &=&Z(X_P)1 
+B^+\left(\sum_{n=0}^{+\infty} np_n X_P^{n-1}
\tilde{Z}(\Delta(X_P))
\right)\\
 &=&Z(X_P)1 
+B^+\left(P'(X_P).
(Z\otimes Id)\circ \Delta(X_P))
\right)
\end{eqnarray*}
We consider the following linear application:
$$ L_P:\left\{\begin{array}{rcl}
\overline{\h} & \longrightarrow & \overline{\h}\\
a& \longrightarrow &B^+(P'(X_P)a).
\end{array}\right. $$
Then, immediately, for all $a \in \overline{\h}$, $val(L_P(a))\geq val(a)+1$,
so  $Id-L_P$ is invertible.
Moreover, by the preceding computation:
\begin{eqnarray*}
&& \tilde{Z} \circ \Delta(X_P)=Z(X_P)1+L_P((Z\otimes Id) \circ \Delta(X_P))\\
&\Longleftrightarrow &
(Id-L_P)((Z\otimes Id) \circ \Delta(X_P))=Z(X_P)1\\
&\Longleftrightarrow &
\tilde{Z} \circ \Delta(X_P)=Z(X_P)(Id-L_P)^{-1}(1).
\end{eqnarray*}
Hence, as $Z(X_P)=1$, lemma \ref{lemme27} induces the following result:

\begin{prop}
\label{prop5}
Let $P \in K[[h]]_1$. If $\A_P$ is a Hopf subalgebra of $\h$, then:
$$(Id-L_P)^{-1}(1) \in \A_P.$$
\end{prop}

\subsection{Proof of $2\Longrightarrow 3$}

We put $\displaystyle Y=\sum_{k=0}^{+\infty} b_n=(Id-L_P)^{-1}(1)$.
Then $b_n$ can be inductively computed in the following way:
\begin{eqnarray}
\label{eq3}
&& \left\{ \begin{array}{rcl}
b_0&=&1,\\
b_{n+1}&=& \displaystyle \sum_{k=1}^{n}\: \sum_{\alpha_1+\ldots+\alpha_k=n} 
(k+1)p_{k+1} B^+(a_{\alpha_1}\ldots a_{\alpha_k})\\
&&\displaystyle +\sum_{k=1}^{n}\: \sum_{\alpha_1+\ldots+\alpha_k=n} 
kp_k B^+(b_{\alpha_1}a_{\alpha_2} \ldots a_{\alpha_k}).
\end{array}
\right.
\end{eqnarray}
In particular, $b_1=p_1 \tun$.\\

Suppose that $\A_P$ is a Hopf subalgebra.
Then $b_n \in (\A_P\cap Vect(\T))_n=Vect(a_n)$ for all $n \geq 1$,
so there exists $\alpha_n \in K$, such that $b_n=\alpha_n a_n$.
Let us compare the projection on $Im((B^+)^2)$ of  $a_{n+1}$ and $b_{n+1}$:
\begin{eqnarray*}
\left\{ \begin{array}{rcl}
p_1B^+(a_n) &\mbox{for}&a_{n+1},\\
2p_2B^+(a_n)+p_1B^+(b_n) =(2p_2+p_1\alpha_n)B^+(a_n)  &\mbox{for}&b_{n+1}.
\end{array}
\right.
\end{eqnarray*}
Suppose that $p_1\neq 0$. Then the $a_n$'s are all non zero by lemma \ref{lemme29}, so
$\alpha_n$ is uniquely determined for all $n \in \mathbb{N}^*$.
We then obtain, by comparing the projections of $a_{n+1}$ and $b_{n+1}$ over $Im((B^+)^2)$:

$$\left\{ \begin{array}{rcl}
\alpha_1&=&p_1,\\
\alpha_{n+1}&=&2\displaystyle \frac{p_2}{p_1}+\alpha_n.
\end{array}
\right.$$
Hence, for all $n \in \mathbb{N}^*$, $\displaystyle \alpha_n=p_1+2\frac{p_2}{p_1}(n-1)$.\\

Let us compare the coefficient of $B^+(\tun^n)$ in $a_{n+1}$ and in $b_{n+1}$ with (\ref{eq1})
and (\ref{eq3}). we obtain:
\begin{eqnarray*}
\left\{ \begin{array}{rcl}
p_n &\mbox{for}&a_{n+1},\\
(n+1)p_{n+1}+np_np_1&\mbox{for}&b_{n+1}.
\end{array}
\right.
\end{eqnarray*}
Hence, $\alpha_{n+1} p_n=(n+1)p_{n+1}+np_np_1$ for all $n \geq 1$.
As a consequence:
$$(n+1)p_{n+1}+ \left(p_1-2\frac{p_2}{p_1}\right)np_n=p_1 p_n.$$
This property is still true for $n=0$, as $p_0=1$.
By multiplying by $h^n$ and taking the sum:
$$P'(h)+\left( p_1-2\frac{p_2}{p_1}\right)hP'(h)=p_1P(h).$$
We then put $\alpha =p_1$ and $\displaystyle \beta= 2\frac{p_2}{p_1^2}-1$.
Hence:
$$(1-\alpha\beta h)P'(h)=\alpha P(h).$$
This equality is still true if $p_1=0$, with $\alpha=0$ and any $\beta$.
Hence, we have shown:
\begin{prop}
\label{prop8}
If $\A_P$ is a Hopf subalgebra of $\h$, then there exists  $(\alpha,\beta)\in K^2$, such that
$P$ satisfies the following differential system:
$$ S_{\alpha,\beta}: \left\{\begin{array}{rcl}
(1-\alpha\beta h)P'(h)&=&\alpha P(h)\\
P(0)&=&1.
\end{array}\right. $$
\end{prop}
This implies:
\begin{enumerate}
\item $P(h)=1$ if $\alpha=0$.
\item $P(h)=e^{\alpha h}$ if $\beta=0$.
\item $\displaystyle P(h)=(1-\alpha\beta h)^{-\frac{1}{\beta}}$ if $\alpha \beta \neq 0$.
\end{enumerate}

\section{Is $\A_{N,\alpha,\beta}$ a Hopf subalgebra?}

The aim of this section is to prove in $4 \Longrightarrow 1$ in theorem \ref{theoprincipal}.

\subsection{Generators of $\A_{N,\alpha,\beta}$.}

We denote by $P_{\alpha,\beta}$ the solution of $S_{\alpha,\beta}$
and we put $\A_{\alpha,\beta}=\A_{P_{\alpha,\beta}}$
and $\A_{N,\alpha,\beta}=\A_{N,P_{\alpha,\beta}}$,
in order to simplify the notations. 
We  also put $\X_{P_{\alpha,\beta}}=\displaystyle \sum_{n\in \mathbb{N}} \a_n(\alpha,\beta)$
and $\displaystyle P_{\alpha,\beta}=\sum_{n=0}^{+\infty} p_n(\alpha,\beta) h^n$.
The system $S_{\alpha,\beta}$ is equivalent to:
$$ S'_{\alpha,\beta}: \left\{\begin{array}{rcl}
p_0(\alpha,\beta)&=&1,\\
p_{n+1}(\alpha,\beta)&=&\displaystyle \alpha \frac{1+n\beta}{n+1}p_n(\alpha,\beta)
\mbox{ for all $n \in \mathbb{N}$}.
\end{array}\right. $$

\begin{defi}
\textnormal{
\begin{enumerate}
\item For all $i\in \mathbb{N}^*$, we put $[i]_\beta=(1+\beta(i-1))$. 
In particular, $[i]_1=i$ et $[i]_0=1$. 
\item We put $[i]_\beta!=[1]_\beta\ldots [i]_\beta$.
In particular, $[i]_1!=i!$ et $[i]_0!=1$. We also put $[0]_\beta!=1$.
\end{enumerate}}
\end{defi}

Immediately, for all $n\in\mathbb{N}$:
$$p_n(\alpha,\beta)=\alpha^n\frac{[n]_\beta!}{n!}.$$

For all $F \in \F$, we define the coefficient $F!$ by:
$$F!=\prod_{\mbox{$s$ vertex of $F$}} (\mbox{fertility of $s$})!.$$
Note that these are not the coefficients $F!$ defined in \cite{Brouder,Hoffman,Stanley}.
They  can be inductively defined by:
$$\left\{\begin{array}{rcl}
 \tun!&=&1,\\
(t_1 \ldots t_k)!&=&t_1!\ldots t_k!,\\
B^+(F)!&=&k! F!
\end{array} \right.$$
In a similar way, we define the following coefficients:
$$[F]_\beta!=\prod_{\mbox{$s$ vertex of $F$}} [\mbox{fertility of $s$}]_\beta!.$$
They can also be inductively defined:
$$\left\{\begin{array}{rcl}
 [\tun]_\beta!&=&1,\\
\left[t_1\ldots t_k\right]_\beta!&=&[t_1]_\beta!\ldots [t_k]_\beta!,\\
\left[B^+(F)\right]_\beta!&=&[k]_\beta! [F]_\beta!
\end{array} \right.$$
In particular, for all forest $F$, $[F]_1!=F!$ and $[F]_0!=1$.\\

Finally, for all $n \in \mathbb{N}^*$, we put
$\displaystyle \a_n(\alpha,\beta)=\sum_{t\in \TT,\: |t|=n} \a_t(\alpha,\beta) t$.

\begin{theo}
\label{theo12}
For any tree $t$,
$$\a_t(\alpha,\beta)=\alpha^{|t|-1}\frac{[t]_\beta!}{t!}.$$
In particular, $\a_t(1,0)=\frac{1}{t!}$, $\a_t(1,1)=1$ and $\a_t(0,\beta)=\delta_{t,\tun}$
for all $\beta$. Moreover:
$$ \a_t(1,-1)=\left\{\begin{array}{cl}
 0& \mbox{ if $t$ is not a ladder,}\\
1 &\mbox{ if $t$ is a ladder.}
\end{array}\right. $$
\end{theo}

{\bf Proof.}
 Induction on $|t|$. If $|t|=1$, then $t=\tun$
and $\a_t(\alpha,\beta)=1$. Suppose the result true for all tree
of weight strictly smaller than $|t|$. 
Then, with $t=B^+(t_1\ldots t_k)$, by (\ref{eq1bis}):
\begin{eqnarray*}
\a_t(\alpha,\beta)&=&
\alpha^{|t_1|-1+\ldots+|t_k|-1} \frac{[t_1]_\beta!\ldots [t_k]_\beta!}
{t_1!\ldots t_k!}\alpha^k\frac{[k]_\beta!}{k!}\\
&=&\alpha^{|t_1|+\ldots+|t_k|}\frac{[t]_\beta!}{t!}\\
&=&\alpha^{|t|-1}\frac{[t]_\beta!}{t!}.
\end{eqnarray*}
The formulas for $(\alpha,\beta)=(1,0)$, $(1,1)$ and $(0,\beta)$ are easily deduced.
Finally, for $(\alpha,\beta)=(1,-1)$, it is enough to observe that
$[1]_\beta!=1$ and $[k]_\beta!=0$ if $k\geq 2$. $\Box$\\

{\bf Examples.}
\begin{eqnarray*}
\a_1(\alpha,\beta) &=&\tun\\
\a_2(\alpha,\beta) &=&\alpha \tdeux \\
\a_3(\alpha,\beta) &=&\alpha^2\left(\frac{(1+\beta)}{2} \ttroisun+\ttroisdeux \right)\\
\a_4(\alpha,\beta) &=&\alpha^3\left(\frac{(1+2\beta)(1+\beta)}{6} \tquatreun
+\frac{(1+\beta)}{2} \tquatredeux+\frac{(1+\beta)}{2} \tquatretrois
+\frac{(1+\beta)}{2}\tquatrequatre+\tquatrecinq\right)\\
\a_5(\alpha,\beta) &=&\alpha^4
\left(
\begin{array}{c}
\frac{(1+3\beta)(1+2\beta)(1+\beta)}{24}\: \tcinqun
+\frac{(1+2\beta)(1+\beta)}{6} \tcinqdeux 
+\frac{(1+2\beta)(1+\beta)}{6} \tcinqtrois 
+\frac{(1+2\beta)(1+\beta)}{6} \tcinqquatre\\
+\frac{(1+\beta)^2}{4}\tcinqsix
+\frac{(1+\beta)^2}{4}\tcinqsept
+\frac{(1+\beta)}{2}\tcinqhuit 
+\frac{(1+\beta)}{2}\tcinqneuf 
+\frac{(1+2\beta)(1+\beta)}{6}\tcinqdix\\
+\frac{(1+\beta)}{2}\tcinqcinq
+\frac{(1+\beta)}{2}\tcinqonze 
+\frac{(1+\beta)}{2}\tcinqdouze 
+\frac{(1+\beta)}{2}\tcinqtreize +\tcinqquatorze
\end{array}\right)
\end{eqnarray*}

In particular, $\a(1,1)$ is the sum of all planar trees of weight $n$,
so $\A_{N,1,1}$ is the subalgebra of formal diffeomorphisms described in \cite{Foissy3}.
Moreover, $\a(1,-1)$ is the ladder of weight $n$, so $\A_{N,1,-1}$ is the subalgebra
of ladders of $\H$.

\subsection{Equalities of the subalgebras $\A_{N,P}$}

\begin{lemme}
Let $P,Q \in K[[h]]_1$. Suppose that $Q(h)=P(\gamma h)$ for a certain $\gamma$.
We denote $\displaystyle \X_P=\sum_{n\geq 1} \a_n$. 
Then $\displaystyle \X_Q=\sum_{n\geq 1} \gamma^{n-1} \a_n$.
In particular, if $\gamma \neq 0$, $\A_{N,P}=\A_{N,Q}$.
\end{lemme}

{\bf Proof.} We put $\displaystyle \Y=\sum_{n\geq 1} \gamma^{n-1} a_n$.
Then:
\begin{eqnarray*}
B^+(Q(\Y)) &=&\sum_{n\in \mathbb{N}}\:\sum_{k=1}^n \:\sum_{n_1+\ldots+n_k=n}\gamma^k p_k
B^+(\gamma^{n_1-1}\a_{n_1}\ldots \gamma^{n_k-1}\a_{n_k})\\
 &=&\sum_{n\in \mathbb{N}}\:\sum_{k=1}^n \:\sum_{n_1+\ldots+n_k=n}\gamma^{k+n-k} p_k
B^+(\a_{n_1}\ldots \a_{n_k})\\
 &=&\sum_{n\in \mathbb{N}}\:\gamma^n\sum_{k=1}^n \:\sum_{n_1+\ldots+n_k=n} p_k
B^+(\a_{n_1}\ldots \a_{n_k})\\
 &=&\sum_{n\in \mathbb{N}}\:\gamma^n \a_{n+1}\\
&=&\Y.
\end{eqnarray*}
By unicity, $\Y=\X_Q$. $\Box$

\begin{theo}
\label{theoegalite}
Let $(\alpha,\beta)$ and $(\alpha',\beta') \in K^2$. The following assertions are equivalent:
\begin{enumerate}
\item $\A_{N,\alpha,\beta}=\A_{N,\alpha',\beta'}$.
\item $\A_{\alpha,\beta}=A_{\alpha',\beta'}$.
\item ($\beta=\beta'$ and $\alpha\alpha'\neq 0$) or $(\alpha=\alpha'=0)$.
\end{enumerate}
\end{theo}

{\bf Proof.}

$1 \Longrightarrow 2$. Obvious.\\

$2 \Longrightarrow 3$. By theorem \ref{theo12}:
$$\left\{ \begin{array}{rclcrclcrcl}
a_1&=&\tun,&& a_2&=&\alpha \tdeux,&&
a_3&=&\alpha^2 \ttroisdeux+\alpha^2\frac{(1+\beta)}{2}\ttroisun;\\ \\
a'_1&=&\tun,&& a'_2&=&\alpha' \tdeux,&&
a'_3&=&\alpha'^2 \ttroisdeux+\alpha'^2\frac{(1+\beta')}{2}\ttroisun.
\end{array}
\right.$$
As $\A_{\alpha,\beta}=\A_{\alpha',\beta'}$,
there exists $\gamma \neq 0$, sucht that $\alpha'\tdeux=\gamma \alpha \tdeux$.
Hence, $\alpha'=\gamma \alpha$. In particular, if $\alpha=0$, then $\alpha'=0$.
Suppose that $\alpha \neq 0$. As $a_3$ and $a'_3$ are colinear, the following
determinant is zero:
$$\begin{array}{|cc|}
\alpha^2& \alpha^2\frac{(1+\beta)}{2}\\
\alpha'^2& \alpha'^2\frac{(1+\beta')}{2}\\
\end{array}
=\frac{1}{2}\alpha^2\alpha'^2(\beta'-\beta)=0.$$
As $\alpha$ and $\alpha'$ are non zero, $\beta=\beta'$. \\

$3 \Longrightarrow 1$. Suppose first $\alpha=\alpha'=0$.
Then $P_{\alpha,\beta}=P_{\alpha',\beta'}=1$, so $\A_{N,\alpha,\beta}=\A_{N,\alpha',\beta'}$.
Suppose $\beta=\beta'$ and $\alpha \alpha' \neq 0$. Then there exists $\gamma\in K-\{0\}$, 
such that $\alpha=\gamma \alpha'$.
Then, immediately, $P_{\alpha,\beta}(\gamma h)=P_{\alpha',\beta'}(h)$.
By the preceding lemma, $\A_{N,\alpha,\beta}=\A_{N,\alpha',\beta'}$. $\Box$

\subsection{The $\A_{N,\alpha,\beta}$'s are Hopf subalgebras}

We now prove $4 \Longrightarrow 1$ in theorem \ref{theoprincipal}. 
If $\alpha=0$, then $\A_{N,\alpha,\beta}=K[\tun]$ and it is obvious.
We take $\alpha \neq 0$. By theorem \ref{theoegalite},
we can suppose that $\alpha=1$.

\begin{lemme}
\label{lemmapoly}
Let $k,n \in \mathbb{N}^*$. We consider the following
element of $K[X_1,\ldots,X_n]$:
$$P_k(X_1,\ldots,X_n)=\sum_{\alpha_1+\ldots+\alpha_n=k}
\frac{X_1(X_1+1)\ldots (X_1+\alpha_1-1)}{\alpha_1!}
\ldots \frac{X_n(X_n+1)\ldots (X_n+\alpha_n-1)}{\alpha_n!}.$$
By putting $S=X_1+\ldots+X_n$:
$$P_k(X_1,\ldots,X_n)=\frac{S(S+1)\ldots(S+k-1)}{k!}.$$
\end{lemme}

{\bf Proof.} Induction on $k$. This is obvious for $k=1$. Suppose the result true at rank $k$.
Then:
\begin{eqnarray*}
&&P_k(X_1,\ldots,X_n)(X_1+\ldots+X_n+k)\\
&=&\sum_{\alpha_1+\ldots+\alpha_n=k}\:
\sum_{i=1}^n \left(
\begin{array}{c}
\frac{X_1(X_1+1)\ldots (X_1+\alpha_1-1)}{\alpha_1!}\\
\vdots\\
\frac{X_i(X_i+1)\ldots (X_i+\alpha_i)}{\alpha_i!}\\
\vdots\\
\frac{X_n(X_n+1)\ldots (X_n+\alpha_n-1)}{\alpha_n!}
\end{array}\right)\\
&=&\sum_{\alpha'_1+\ldots+\alpha'_n=k+1}\:
\sum_{i=1}^n \alpha'_i \left(
\begin{array}{c}
\frac{X_1(X_1+1)\ldots (X_1+\alpha'_1-1)}{\alpha'_1!}\\
\vdots\\
\frac{X_i(X_i+1)\ldots (X_i+\alpha'_i-1)}{\alpha'_i!}\\
\vdots\\
\frac{X_n(X_n+1)\ldots (X_n+\alpha'_n-1)}{\alpha'_n!}
\end{array}\right)\\
&=&(k+1)\sum_{\alpha'_1+\ldots+\alpha'_n=k+1}\:
\left(
\begin{array}{c}
\frac{X_1(X_1+1)\ldots (X_1+\alpha'_1-1)}{\alpha'_1!}\\
\vdots\\
\frac{X_i(X_i+1)\ldots (X_i+\alpha'_i-1)}{\alpha'_i!}\\
\vdots\\
\frac{X_n(X_n+1)\ldots (X_n+\alpha'_n-1)}{\alpha'_n!}
\end{array}\right)\\
&=&(k+1)P_{k+1}(X_1,\ldots,X_n).
\end{eqnarray*}
This implies the announced result. $\Box$\\

Let $F=t_1\ldots t_k$ be a forest and $t$ be a tree. 
Using the dual basis $(e_F)_{F\in \FF}$:
\begin{eqnarray*}
\mbox{coefficient of $F \otimes t$ in $\Delta(\X_{1,\beta})$}
&=&<e_F \otimes e_t,\Delta(\X_{1,\beta})>\\
&=&<e_Fe_t,\X_{1,\beta}>\\
&=&\sum_{\mbox{$s$ grafting of $F$ on $t$}} <e_s,\X_{1,\beta}>\\
&=&\sum_{\mbox{$s$ tree, grafting of $F$ on $t$}}  \frac{[s]_\beta!}{s!}.
\end{eqnarray*}
Let $n$ be the weight of $t$ and $s_1,\ldots,s_n$ its vertices.
Let $f_i$ be the fertility of $s_i$. Let $(\alpha_1,\ldots,\alpha_n)$
such that $\alpha_1+\ldots+\alpha_n=k$ and consider the graftings 
of $F$ on $t$ such that $\alpha_i$ trees of $F$ are grafted on $s_i$ for all $i$. Then:
\begin{enumerate}
\item If $s$ is such a grafting, we have:
$$  \left\{ \begin{array}{rcl}
[s]_\beta!&=&\displaystyle [t]_\beta![t_1]_\beta!\ldots[t_k]_\beta!
\frac{[f_1+\alpha_1]_\beta!}{[f_1]_\beta!}\ldots \frac{[f_n+\alpha_n]_\beta!}{[f_n]_\beta!},\\[4mm]
s!&=&\displaystyle t!t_1!\ldots t_k! \frac{(f_1+\alpha_1)!}{f_1!}\ldots \frac{(f_n+\alpha_n)!}{f_n!}.
\end{array}
\right.$$
\item The number of such graftings is:
$$\binom{f_1+\alpha_1}{\alpha_1}
\ldots \binom{f_n+\alpha_n}{\alpha_n}.$$
\end{enumerate}
Hence, by putting $x_i=f_i+1/\beta$ and $s=x_1+\ldots+x_n$,
by lemma \ref{lemmapoly}:
\begin{eqnarray*}
&&\mbox{coefficient of $F \otimes t$ in $\Delta(\X_{1,\beta})$}\\
&=&\sum_{\alpha_1+\ldots+\alpha_n=k} \frac{[t]_\beta!}{t!}
\frac{[t_1]_\beta!}{t_1!}\ldots \frac{[t_k]_\beta!}{t_k!}
\frac{[f_1+\alpha_1]_\beta!}{[f_1]_\beta!\alpha_1!}
\ldots \frac{[f_n+\alpha_n]_\beta!}{[f_n]_\beta!\alpha_n!}\\
&=&\frac{[t]_\beta!}{t!} \frac{[t_1]_\beta!}{t_1!}\ldots \frac{[t_k]_\beta!}{t_k!}
\sum_{\alpha_1+\ldots+\alpha_n=k}\: \prod_{i=1}^n
\frac{(1+f_i\beta)\ldots(1+(f_i+\alpha_i-1)\beta)}{\alpha_i!}\\
&=&\frac{[t]_\beta!}{t!} \frac{[t_1]_\beta!}{t_1!}\ldots \frac{[t_k]_\beta!}{t_k!}
\sum_{\alpha_1+\ldots+\alpha_n=k}\: \prod_{i=1}^n
\beta^{\alpha_i} \frac{x_i(x_i+1)\ldots (x_i+\alpha_i-1)}{\alpha_i}\\
&=&\frac{[t]_\beta!}{t!} \frac{[t_1]_\beta!}{t_1!}\ldots \frac{[t_k]_\beta!}{t_k!}
\beta^k P_k(x_1,\ldots,x_n)\\
&=&\frac{[t]_\beta!}{t!} \frac{[t_1]_\beta!}{t_1!}\ldots \frac{[t_k]_\beta!}{t_k!}
\beta^k \frac{s(s+1)\ldots(s+k-1)}{k!}.
\end{eqnarray*}
Moreover, as $t$ is a tree:
$$s=f_1+\ldots+f_n+n/\beta
=\mbox{number of edges of $t$}+n/\beta
=n-1+n/\beta
=n(1+1/\beta)-1.$$
So, as $Q_k(S)=\displaystyle \frac{S(S+1)\ldots (S+k-1)}{k!}$:
\begin{eqnarray*}
\Delta(\X_{1,\beta})&=&\X_{1,\beta} \otimes 1+\sum_{k=0}^{\infty} \sum_{F=t_1\ldots t_k,\:t}
\frac{[t]_\beta!}{t!} \frac{[t_1]_\beta!}{t_1!}\ldots \frac{[t_k]_\beta!}{t_k!}
\beta^k Q_k(|t|(1+1/\beta)-1) F \otimes t\\
&=&\X_{1,\beta} \otimes 1+\sum_{n=1}^{\infty}\sum_{k=0}^{\infty} 
Q_k(n(1+1/\beta)-1)\beta^k \X_{1,\beta}^k\otimes \a_n(1,\beta)\\
&=&\X_{1,\beta} \otimes 1+\sum_{n=1}^{\infty}
(1-\beta \X_{1,\beta})^{-n(1/\beta+1)+1 }\otimes \a_n(1,\beta).
\end{eqnarray*}

\begin{prop}
\label{propcoproduit}
The coproduct of the $\a_n(1,\beta)$'s is given by:
$$\Delta(\X_{1,\beta})=\X_{1,\beta} \otimes 1+\sum_{n=1}^{\infty}
(1-\beta \X_{1,\beta})^{-n(1/\beta+1)+1 }\otimes \a_n(1,\beta).$$
As a consequence, $\A_{N,1,\beta}$ is a Hopf subalgebra of $\H$.
\end{prop}

{\bf Remark.} By taking the abelianization of $\A_{N,1,\beta}$, the same
holds in $\A_{1,\beta}$:
$$\Delta(X_{1,\beta})=X_{1,\beta} \otimes 1+\sum_{n=1}^{\infty}
(1-\beta X_{1,\beta})^{-n(1/\beta+1)+1 }\otimes a_n(1,\beta).$$

\section{Isomorphisms between the $\A_{N,\alpha,\beta}$'s}

\subsection{Another system of generators of $\A_{N,1,\beta}$}

{\bf Notation.} We denote by $B^-$ the inverse of $B^+:\H\longrightarrow Vect(\TT)$,
that is to say the application defined on a tree by deleting the root.\\

We define $\b_n(\alpha,\beta)=B^-(\a_{n+1}(\alpha,\beta))$ for all $n \in \mathbb{N}$, and:
$$\Y(\alpha,\beta)=\sum_{n=0}^\infty \b(\alpha,\beta).$$
We have:
\begin{eqnarray*}
\Y(\alpha,\beta)&=&\sum_{F \in \F}
\alpha^{|B^+(F)|-1} \frac{[B^+(F)]_\beta!}{B^+(F)!} F\\
&=&\sum_{k=0}^\infty \sum_{t_1,\ldots,t_k \in \T} \alpha^{|t_1|+\ldots+|t_k|}
\frac{[k]_\beta! [t_1]_\beta! \ldots [t_k]_\beta!}{k! t_1!\ldots t_k!} t_1\ldots t_k\\
&=&\sum_{k=0} \frac{[k]_\beta!}{k!} \X_{\alpha,\beta}^k\\
&=&\sum_{k=0}^\infty \frac{1(1+\beta)\ldots(1+\beta(k-1))}{k!} \X_{\alpha,\beta}^k\\
&=&\sum_{k=0}^\infty \frac{1/\beta(1/\beta+1)\ldots(1/\beta+k-1))}{k!} \beta^k \X_{\alpha,\beta}^k\\
&=&\sum_{k=0}^\infty Q_k(1/\beta)\beta^k \X_{\alpha,\beta}^k\\
&=&\left(1-\beta \X_{\alpha,\beta}\right)^{-1/\beta}.
\end{eqnarray*}
By the last equality, $\b_n(\alpha,\beta) \in \A_{N,\alpha,\beta}$ for all $n \in \mathbb{N}$.
Moreover, by the second equality, if $n \geq 1$:
\begin{eqnarray*}
\b_n(\alpha,\beta)&=&\sum_{t\in \T,\:|t|=n} \alpha^n \frac{[t]_\beta!}{t!}t 
+\mbox{forests with more than two trees}\\
&=&\alpha \a_n(\alpha,\beta)+
\mbox{forests with more than two trees}.
\end{eqnarray*}
So $(\b_n(\alpha,\beta))_{n\geq 1}$ is a set of generators of $\A_{N,\alpha,\beta}$
if $\alpha \neq 0$.

\begin{prop}
Suppose $\alpha=1$. Then:
$$\Delta(\Y(1,\beta))=\sum_{n=0}^\infty \Y(1,\beta)^{n(\beta+1)+1}\otimes \b_n(1,\beta).$$
\end{prop}

{\bf Proof.} As $\Y(1,\beta)=B^-(\X(\alpha,\beta))$, by (\ref{E1}):
\begin{eqnarray*}
\Delta(\Y(1,\beta))&=&(Id \otimes B^-)\left(\Delta(\X(1,\beta))-\X(1,\beta)\otimes 1\right)\\
&=&\sum_{n=1}^{\infty}(1-\beta \X_{1,\beta})^{-n(1/\beta+1)+1 }\otimes B^-(\a_n(1,\beta))\\
&=&\sum_{n=1}^{\infty}\Y_{1,\beta}^{n(1+\beta)-\beta }\otimes \b_{n-1}(1,\beta)\\
&=&\sum_{n=0}^{\infty}\Y_{1,\beta}^{(n+1)(1+\beta)-\beta }\otimes \b_n(1,\beta)\\
&=&\sum_{n=0}^{\infty}\Y_{1,\beta}^{n(1+\beta)+1}\otimes \b_n(1,\beta). \:\Box
\end{eqnarray*}

\subsection{Isomorphisms between the $\A_{N,\alpha,\beta}$'s}

\begin{prop}
If $\beta\neq -1$ and $\beta'\neq -1$, then $\A_{N,1,\beta}$
and $\A_{N,1,\beta'}$ are isomorphic.
\end{prop}

{\bf Proof.} Let $\gamma \in K-\{0\}$. We put:
$$\Z(1,\beta)=\Y(1,\beta)^\gamma=\sum_{k=0}^\infty \c_n(1,\beta),$$
with $\c_n(1,\beta) \in \A(1,\beta)$, homogeneous of degree $n$.
This makes sense, because $\b_0(1,\beta)=1$.
Moreover, for all $n \geq 1$:
\begin{eqnarray*}
\c_n(1,\beta)&=&Q_1(\gamma) \b_n(1,\gamma)
+\mbox{forests with more than two trees}\\
&=&\gamma \b_n(1,\gamma)
+\mbox{forests with more than two trees},
\end{eqnarray*}
so $(\c_n(1,\beta))_{n \geq 1}$ is a set of generators of $\A(1,\beta)$.
Moreover:
\begin{eqnarray*}
\Delta(\Z(1,\beta))&=&\sum_{k=0}^\infty Q_k(\gamma)
\left(\sum_{n=1}^\infty \Y(1,\beta)^{n(\beta+1)+1}\otimes \b_n(1,\beta) \right)^k\\
&=&\sum_{k=0}^\infty Q_k(\gamma)
\sum_{a_1,\ldots,a_k\geq 1}\Y(1,\beta)^{(a_1+\ldots+a_k)(\beta+1)+k}
\otimes \b_{a_1}(1,\beta)\ldots \b_{a_k}(1,\beta) \\
&=&\sum_{l=0}^\infty \sum_{k=0}^\infty Q_k(\gamma)
\sum_{a_1+\ldots+a_k=l}\Y(1,\beta)^{l(\beta+1)+k}
\otimes \b_{a_1}(1,\beta)\ldots \b_{a_k}(1,\beta) \\
&=&\sum_{l=0}^\infty \Y(1,\beta)^{l(\beta+1)}
\left(\sum_{k=0}^\infty \sum_{a_1+\ldots+a_k=l} Q_k(\gamma)
\Y(1,\beta)\otimes \b_{a_1}(1,\beta)\ldots \Y(1,\beta)\otimes \b_{a_k}(1,\beta)\right)\\
&=&\sum_{l=0}^\infty \Y(1,\beta)^{l(\beta+1)} \Y(1,\beta)^\gamma \otimes \c_l(1,\beta)\\
&=&\sum_{l=0}^\infty \Z(1,\beta)^{l\left(\frac{\beta+1}{\gamma}\right)+1}\otimes \c_l(1,\beta).
\end{eqnarray*}

We now chose $\gamma=\frac{\beta+1}{\beta'+1}$. As $\beta'\neq -1$, this is well defined;
as $\beta\neq -1$, this is non zero. Then:
$$\Delta(\Z(1,\beta))=
\sum_{l=0}^\infty \Z(1,\beta)^{l\left(\beta'+1\right)+1}\otimes \c_l(1,\beta).$$
So the unique isomorphism of algebras defined by:
$$ \left\{\begin{array}{rcl}
\A_{1,\beta'} & \longrightarrow & \A_{1,\beta}\\
\b_n(1,\beta')  & \longrightarrow & \c_n(1,\beta)
\end{array}\right. $$
is a Hopf algebra isomorphism. $\Box$ \\

In the non commutative case, the following result holds:

\begin{cor}
There are three isomorphism classes of $\A_{N,\alpha,\beta}$'s:
\begin{enumerate}
\item the $\A_{N,1,\beta}$'s, with $\beta\neq -1$.
These are not commutative and not cocommutative.
\item $\A_{N,1,-1}$, isomorphic to $\bf QSym$,
the Hopf algebra of quasi-symmetric functions (\cite{Malvenuto,Stanley2})
This one is not commutative and cocommutative.
\item $\A_{N,0,1}=K[\tun]$. This one is commutative and cocommutative.
\end{enumerate}
\end{cor}

Consequently, in the commutative case:
\begin{cor}
There are three isomorphism classes of $\A_{\alpha,\beta}$'s:
\begin{enumerate}
\item the $\A_{1,\beta}$'s, with $\beta\neq -1$.
These are isomorphic to the Fa\`a di Bruno algebra on one variable. 
\item $\A_{1,-1}$, isomorphic to $Sym$,
the Hopf algebra of symmetric functions.
This one is commutative and cocommutative.
\item $\A_{0,1}=K[\tun]$. 
\end{enumerate}
\end{cor}

{\bf Proof.} As $\A_{\alpha,\beta}$ is the abelianization of $\A_{N,\alpha,\beta}$,
if $\A_{N,\alpha,\beta}\approx \A_{N,\alpha',\beta'}$, then
$\A_{\alpha,\beta}\approx \A_{\alpha',\beta'}$. Moreover,
$\A_{1,\beta}$ is not cocommutative if $\beta\neq -1$, whereas
$\A_{1,-1}$ is. So $A_{1,\beta}$ and $A_{1,-1}$ are not isomorphic if $\beta \neq -1$.
It remains to show that $\A_{1,\beta}$ is isomorphic to the Fa\`a di Bruno Hopf algebra
on one variable if $\beta \neq -1$.
Let us consider the dual Hopf algebra of $\A_{1,\beta}$.
By Cartier-Quillen-Milnor-Moore's theorem (\cite{Milnor}),
this is an enveloping algebra ${\cal U}(\L_{1,\beta})$. 
Moreover, $\L_{1,\beta}$ has for basis $(T_n)_{n\in \mathbb{N}^*}$ defined by:
$$ T_n : \left\{ \begin{array}{rcl}
\L_{1,\beta} &\longrightarrow &  K\\
a_1^{\alpha_1} \ldots a_k^{\alpha_k} &\longrightarrow & 0 \mbox{ if }\alpha_1+\ldots+\alpha_k \neq 1,\\
a_m&\longrightarrow & \delta_{m,n}.
\end{array}
\right.$$
Moreover, $T_n$ is homogeneous of degree $n$.
By proposition \ref{propcoproduit}, for all $i,j \geq 1$:
\begin{eqnarray*}
(T_i\otimes T_j) \circ \Delta(X_{1,\beta})
&=&T_i(X_{1,\beta}) T_j(1)+\sum_{n=1}^{\infty}
\sum_{k=0}^\infty Q_k(n(1/\beta+1)-1)\beta^k T_i(X_{1,\beta}^k) T_j(a_n(1,\beta))\\
&=&0+Q_1(j(1/\beta+1)-1)\beta T_i(X_{1,\beta})\\
&=&j(1+\beta)-\beta.
\end{eqnarray*}
By homogeneity, there exists $\lambda_{i,j} \in K$ such that 
$[T_i,T_j]=\lambda_{i,j}T_{i+j}$. Then:
\begin{eqnarray*}
\lambda_{i,j}&=&[T_i,T_j](X_{1,\beta})\\
&=&(T_i \otimes T_j)\circ  \Delta(X_{1,\beta})-(T_j \otimes T_i)\circ  \Delta(X_{1,\beta})\\
&=&j(1+\beta)-\beta -i(1+\beta)+\beta\\
&=&(i-j)(1+\beta).
\end{eqnarray*}
Then, there exists a Lie algebra morphism:
$$\left\{\begin{array}{rcl}
\L_{\alpha,\beta} & \longrightarrow & Prim(\hfdb ^*)\\
T_n& \longrightarrow &(1+\beta)Z_n.
\end{array}\right. $$
In particular, if $\beta \neq -1$, this is an isomorphism. 
Hence, $A_{1,\beta}^*$ is isomorphic to $\hfdb^*$, so
$A_{1,\beta}$ is isomorphic to the Fa\`a di Bruno Hopf algebra on one variable. $\Box$\\
 
{\bf Remark.} The Connes-Moscovici subalgebra ${\cal H}_{CM}$ of $\h$ (see \cite{Connes,Moscovici})
does not appear here: as it is generated by $\tun, \tdeux,\ttroisun+\ttroisdeux,\ldots$,
it would be $A_{(1,1)}$. The fourth generator of $A_{(1,1)}$ is:
$$\tquatreun+2\tquatredeux+\tquatrequatre+\tquatrecinq,$$
whereas the fourth generator of ${\cal H}_{CM}$ is:
$$\tquatreun+3\tquatredeux+\tquatrequatre+\tquatrecinq.$$
So they are different.

\section{The case of the free Fa\`a di Bruno algebra with $D$ variables}

We here fix an integer $D\geq 1$.
We denote by $W$ the set of non empty words in letters $\{1,\ldots,D\}$.

\subsection{Construction}

We now recall the construction of the free Fa\`a di Bruno algebra in $D$ variables (see \cite{Effros}).
Consider the ring of non commutative formal series $K\langle \langle h_1,\ldots,h_D \rangle \rangle$ on $D$ variables.
We consider:
$$G_D=\left\{
\left( \sum_{w\in W} a^{(i)}_{w} h^w\right)_{1\leq i\leq D}\:/\: a^{(i)}_{j}=\delta_{i,j}\right\}.$$
We use the following convention: if $u_1\ldots u_k \in W$, then $h^w=h_{u_1}\ldots h_{u_k}$.
In other terms, $G_D$ is the set of formal diffeomorphisms on $K^D$ which are tangent to the identity
at the origin. This is a group for the composition of formal series. 

Then $\hfdbD$ is the Hopf algebra of functions on the opposite of the group $G_D$. 
Hence, it is the polynomial ring in variable $Y^i_w$, $1\leq i\leq D$, with the convention that
if $w$ has only one letter $j$, then $Y^i_j=\delta_{i,j}$. The coproduct is given in the following
way: for all $f \in \hfdbD$, for all $P,Q \in G_D$,
$$\Delta(f)(P\otimes Q)=f(Q\circ P).$$
In particular, if:
$$P=\left( \sum_{w\in W} a^{(i)}_{w} h^w\right)_{1\leq i\leq D}\mbox{ and }
Q=\left( \sum_{w\in W} b^{(i)}_{w} h^w\right)_{1\leq i\leq D},$$
then:
\begin{eqnarray*}
\Delta(Y^i_w)(P\otimes Q)&=&Y^i_w(Q\circ P)\\
&=&Y^i_w\left(\left(\sum_{u_1\ldots u_k \in W}
b^j_{u_1\ldots u_k} \sum_{w_1,\ldots,w_k \in W} a_{w_1}^{u_1}\ldots a_{w_k}^{u_k} h^{w_1\ldots w_k}
\right)_{1\leq j \leq D}\right)\\
&=&\sum_{u_1\ldots u_k \in W}\sum_{w_1\ldots w_k=w}b^j_{u_1\ldots u_k} a_{w_1}^{u_1}\ldots a_{w_k}^{u_k}.
\end{eqnarray*}
Hence:
$$\Delta\left( Y_w^i\right)
=\sum_{k=1}^n \:\sum_{1\leq u_i \leq D} \:
\sum_{ \begin{array}{c} \\[-7mm]_{w_1,\ldots,w_k \in W},\\[-1mm] _{w_1\ldots w_k=w}\end{array}}
Y^{u_1}_{w_1}\ldots Y^{u_k}_{w_k}\otimes Y^i_{u_1\ldots u_k}.$$
For $D=1$, we recover $\hfdb$.

\subsection{Subalgebras of $\H^\D$}

We now put $\D=\{1,\ldots,D\}^3$. The elements of $\D$ will be denoted in the following
way: $i,(u_1,u_2)$. In the same way, it is possible to construct
a commutative Hopf algebra $\h^\D$ of rooted trees decorated by $\D$,
and a non commutative Hopf algebra $\H^\D$ of planar rooted trees decorated by $\D$.
In both cases, we define, for all $i,(u_1,u_2) \in \D$, a linear endomorphism
$B_{i,(u_1,u_2)}^+$, which sends forest $F$ on the tree obtained by grafting all the trees of $F$ on a common
root decorated by $i,(u_1,u_2)$.

\begin{defi}
\textnormal{
Let $i\in \{1,\ldots,D\}$ and $w=u_1\ldots u_n \in W$. We define an element $Y_w^i \in \H^\D$
inductively on $n$ in the following way:
$$  \left\{\begin{array}{rcl}
 \Y_w^i &=& \delta_{i,w} \mbox{ if $n=1$},\\
 \Y_w^i&=&\displaystyle \sum_{1\leq \alpha,\beta\leq D}\: \sum_{\substack{w_1,w_2\in W,\\w_1w_2=w}}
 B^+_{i,(\alpha,\beta)}\left(Y_{w_1}^\alpha Y_{w_2}^\beta \right)
 \mbox{ if $n\geq 2$}.
 \end{array}\right. $$}
\end{defi}

{\bf Examples.} For $i$, $u_1$, $u_2$, $u_3$ and $u_4$ elements of $\{1,\ldots, D\}$:
\begin{eqnarray*}
Y_{u_1}^i&=&\delta_{i,u_1},\\
Y_{u_1u_2}^i&=&\tdun{$i,(u_1,u_2)$}\hspace{1cm},\\
Y_{u_1u_2u_3}^i&=&\sum_{1\leq \alpha \leq D}\left(\tddeux{$i,(\alpha,u_3)$}{$\alpha,(u_1,u_2)$}
\hspace{1.2cm}+\tddeux{$i,(u_1,\alpha)$}{$\alpha,(u_2,u_3)$}\hspace{1cm}\right),\\
Y_{u_1u_2u_3u_4}^i&=&\sum_{1\leq \alpha,\beta \leq D}\left(
\hspace{1.1cm} \tdtroisun{$i,(\alpha,\beta)$}{$\beta,(u_3,u_4)$}{ \hspace{-1.2cm}$\alpha,(u_1,u_2)$}
\hspace{1.2cm}+\tdtroisdeux{$i,(\alpha,u_4)$}{$\alpha,(u_1,\beta)$}{$\beta,(u_2,u_3)$}
\hspace{1.2cm}+\tdtroisdeux{$i,(\alpha,u_4)$}{$\alpha,(\beta,u_3)$}{$\beta,(u_1,u_2)$}
\hspace{1.2cm}+\tdtroisdeux{$i,(u_1,\alpha)$}{$\alpha,(u_2,\beta)$}{$\beta,(u_3,u_4)$}
\hspace{1.2cm}+\tdtroisdeux{$i,(u_1,\alpha)$}{$\alpha,(\beta,u_4)$}{$\beta,(u_2,u_3)$}\hspace{1cm}\right).
\end{eqnarray*}
An easy induction shows that $Y^i_{u_1\ldots u_n}$ is homogeneous of degree $n-1$.

\begin{theo}
For all $i \in \{1,\ldots,D\}$, $w \in W$ of length $n$:
$$\Delta\left( Y_w^i\right)
=\sum_{k=1}^n \:\sum_{1\leq \alpha_i \leq D} \:
\sum_{\substack{w_1,\ldots,w_k \in W,\\ w_1\ldots w_k=w}}
Y^{\alpha_1}_{w_1}\ldots Y^{\alpha_k}_{w_k}\otimes Y^i_{\alpha_1\ldots \alpha_k}.$$
\end{theo}

{\bf Proof.} By induction on $n$. It is obvious if $n=1$ or $2$. 
Suppose it is true for all rank $<n$. Then:
\begin{eqnarray*}
\Delta(Y^i_w)&=&\sum_{ \alpha,\beta}\: \sum_{w_1w_2=w}
\Delta \circ B^+_{i,(\alpha,\beta)}\left(Y_{w_1}^\alpha Y_{w_2}^\beta \right)\\
&=&Y^i_w\otimes 1
+\sum_{ \alpha,\beta}\: \sum_{w_1w_2=w}\:
\sum_{w_{1,1}\ldots w_{1,k}=w_1} \\
&&\hspace{1cm}\sum_{\alpha_1,\ldots,\alpha_k}\:
\sum_{w_{2,1}\ldots w_{2,l}=w_2} \:
\sum_{\beta_1,\ldots,\beta_l}\:
Y^{\alpha^1}_{w_{1,1}}\ldots Y^{\alpha^k}_{w_{1,k}}
Y^{\beta^1}_{w_{2,1}}\ldots Y^{\beta^l}_{w_{2,l}}
\otimes B^+_{i,(\alpha,\beta)}
\left(Y^\alpha_{\alpha_1\ldots\alpha_k}Y^\beta_{\beta_1\ldots\beta_l}\right)\\
&=&Y^i_w\otimes 1+\sum_{k\geq 2} \: \sum_{w_1\ldots w_k=w} \:\sum_{\alpha_1,\ldots,\alpha_k}
Y_{w_1}^{\alpha_1}\ldots Y_{w_k}^{\alpha_k}\otimes
\left(\sum_{\alpha,\beta}\:\sum_{w'_1w'_2=\alpha_1\ldots \alpha_k}\:
B_{i,(\alpha,\beta)}^+\left(Y_{w'_1}^\alpha Y_{w'_2}^\beta\right) \right)\\
&=&Y^i_w\otimes 1+\sum_{k\geq 2} \: \sum_{w_1\ldots w_k=w} \: \sum_{\alpha_1,\ldots,\alpha_k}
Y_{w_1}^{\alpha_1}\ldots Y_{w_k}^{\alpha_k}\otimes Y_{\alpha_1\ldots \alpha_k}^i\\
&=&\sum_{k\geq 1} \: \sum_{w_1\ldots w_k=w} \: \sum_{\alpha_1,\ldots,\alpha_k}
Y_{w_1}^{\alpha_1}\ldots Y_{w_k}^{\alpha_k}\otimes Y_{\alpha_1\ldots \alpha_k}^i. \: \Box \\
\end{eqnarray*}

Hence, the subalgebra of $\H^\D$ generated by the $Y_w^i$'s is a Hopf subalgebra.
Its abelianization can be seen as a subalgebra of $\h^\D$, and is isomorphic to $\hfdbD$.\\

{\bf Remark.} In the case where $D=1$, we put $Y^1_{\underbrace{1\ldots 1}_{n+1\:times}}=Y_n$.
Then, by definition:
$$  \left\{\begin{array}{rcl}
 Y_0 &=&1,\\
 Y_1&=&\tdun,\\
 Y_n&=&\displaystyle \sum_{k=0}^{n-1} B^+(Y_k Y_{n-1-k})
=2B^+(Y_{n-1})+\sum_{k=1}^{n-2} B^+(\Y_kY_{n-1-k})\mbox{ if }n\geq2.
 \end{array}\right. $$
Hence, by (\ref{eq2}), this is the subalgebra associated to $1+2h+h^2=(1+h)^2=P_{4,-\frac{1}{2}}(h)$,
namely $\A_{D,4,-\frac{1}{2}}$.

\subsection{Description of the $Y_w^i$'s in the generic case}

\begin{defi}
\textnormal{A tree $t \in \T^\D$ is {\it admissible} if:
\begin{enumerate}
\item Every vertex is of fertility less than $2$.
\item For each vertex of fertility $1$, the decorations are set in this way:
$$\tddeux{$i,(a,b) $}{$ a,(c,d) $}
\hspace{.8cm}\mbox{ or } \tddeux{$i,(a,b) $}{$ b,(c,d) $}\hspace{.7cm}.$$
\item For each vertex of fertility $2$, the decorations are set in this way:
$$\tdtroisun{$i,(a,b)$}{$b,(e,f)$}{\hspace{-.8cm}$a,(c,d)$}\hspace{.7cm}.$$
\end{enumerate}}
\end{defi}

Let $t$ be an admissible tree. 
We associate to it a word in $W$ in the following inductive way:
\begin{enumerate}
\item $w(\tdun{$i,(a,b)$} \hspace{.6cm})=ab$.
\item If the root of $t$ has fertility $1$, with decorations set as 
$\tddeux{$i,(a,b) $}{$ a,(c,d) $}\hspace{.7cm}$, then, if we denote $t'=B^-(t)$, $w(t)=w(t')b$.
\item If the root of $t$ has fertility $1$, with decorations set as 
$\tddeux{$i,(a,b) $}{$ b,(c,d) $}\hspace{.7cm}$, then, if we denote $t'=B^-(t)$, $w(t)=aw(t')$.
\item If the root of $t$ has fertility $2$, then, if we denote $t't''=B^-(t)$, $w(t)=w(t')w(t'')$.
\end{enumerate}

{\bf Remark.} The cases $1$ and $2$ are not incompatible, so $w(t)$ is not well defined.
For example, for $t=\tddeux{$i,(a,a) $}{$ a,(c,d) $}\hspace{.7cm}$, two results are possible:
$cda$ and $acd$.\\

An easy induction shows that:

\begin{prop}
Suppose that $w \in W$ is generic, that is to say all his letters are distinct. 
Then $Y_w^i$ is the sum of admissible trees $t$ such that:
\begin{enumerate}
\item $w(t)=w$.
\item The decoration of the root of $t$ is of the form $i,(a,b)$, with $1 \leq ,a,b \leq D$.
\end{enumerate}
\end{prop}

If the word is not generic, we obtain $Y_w^i$ by specializing the generic case. For example,
if $w=aaa$, we have:
\begin{eqnarray*}
Y_{abc}^i&=&\sum_{1\leq \alpha \leq D}\tddeux{$i,(\alpha,c)$}{$\alpha,(a,b)$}
\hspace{1.2cm}+\sum_{1\leq \alpha \leq D}\tddeux{$i,(a,\alpha)$}{$\alpha,(b,c)$}\\
\Longrightarrow \: Y_{aaa}^i&=&\sum_{1\leq \alpha \leq D}\tddeux{$i,(\alpha,a)$}{$\alpha,(a,a)$}
\hspace{1.2cm}+\sum_{1\leq \alpha \leq D}\tddeux{$i,(a,\alpha)$}{$\alpha,(a,a)$}\hspace{.7cm}.
\end{eqnarray*}
In particular, $\tddeux{$i,(a,a)$}{$a,(a,a)$} \hspace{.7cm}$ appears with multiplicity $2$.

\bibliographystyle{amsplain}
\bibliography{biblio}

\end{document}